\documentclass{article}
\usepackage{graphicx}
\usepackage{amsmath}

\input epsf.tex

\newtheorem{theorem}{Theorem}

\newtheorem{lemma}[theorem]{Lemma}

\newtheorem{proposition}{Proposition}

\newenvironment{proof}[1][Proof]{\textbf{#1.} }{\ \rule{0.5em}{0.5em}}

\begin{document}

\title{Extreme coefficients of Jones polynomials and graph theory}
\author{P. M. G. Manch\'on \\
The University of Liverpool}
\maketitle

\begin{abstract}
We find families of prime knot diagrams with arbitrary extreme
coefficients in their Jones polynomials. Some graph theory is
presented in connection with this problem, generalizing ideas by
Yongju Bae and Morton \cite{main} and giving a positive answer to
a question in their paper.
\end{abstract}

\section{Introduction.}

Let $L$ be an oriented link, and $V_L(t)$ its Jones polynomial
with normalization one. We are interested in exhibit examples of
links with arbitrary extreme coefficients in their Jones
polynomials. Consider an unoriented diagram $D$ of $L$. We denote
by $\langle D\rangle$ its Kauffman bracket with normalization
$\langle \raisebox{-.8mm}{\epsfysize.15in \epsffile{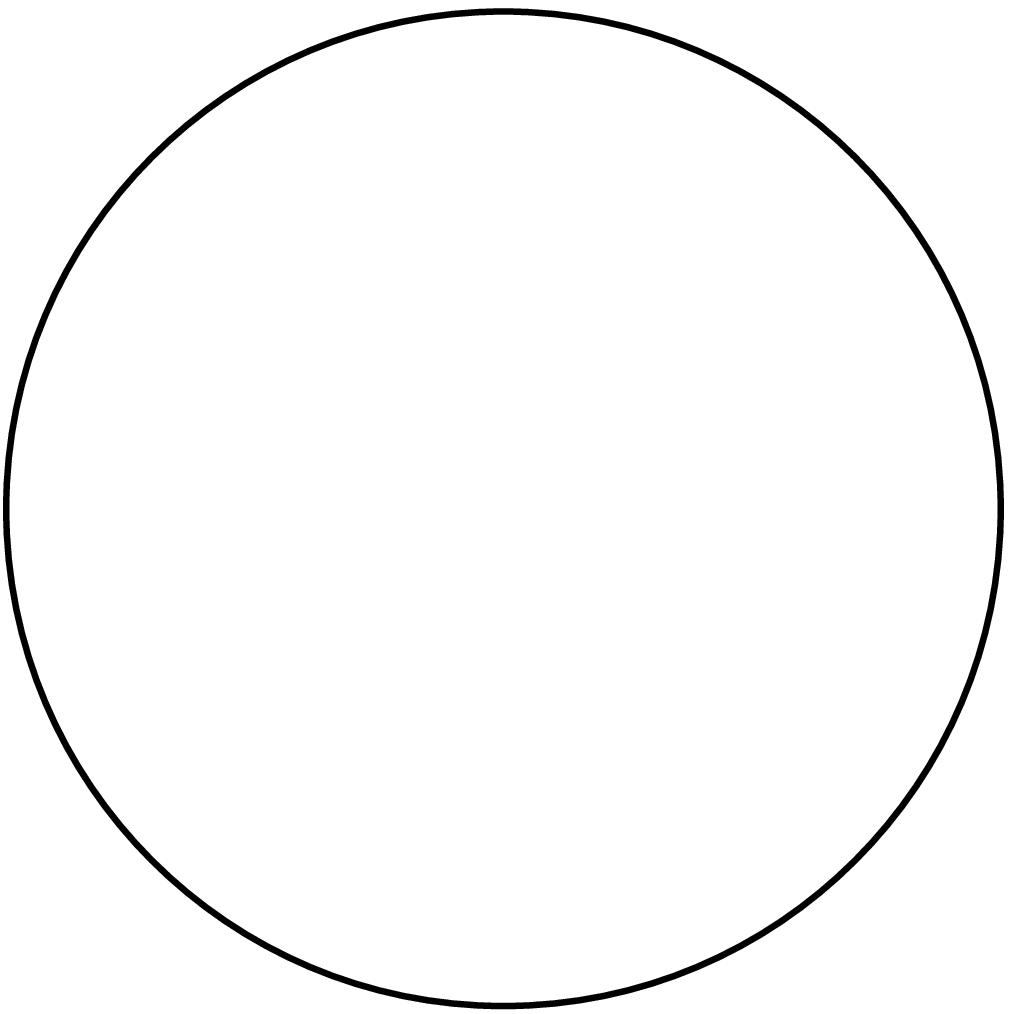}}
\rangle =1$ (see \cite{Lickorish}). Since
$V_L(t)=(-A)^{-3w(D)}\langle D\rangle $ after the substitution
$A=t^{-1/4}$, we have that $span(V_L(t))=span(\langle D\rangle)/4
$ and the coefficients of both polynomials are the same, maybe up
to sign.


Hence we will work with the Kauffman bracket of unoriented
diagrams. We recall the definition of this polynomial using the
states sum:

$$\langle D \rangle =\sum _s\langle D,s \rangle $$

\noindent where the sum is taken over all states $s$ of $D$ and
$\langle D,s \rangle =A^{a(s)-b(s)}(-A^{-2}-A^2)^{|s|-1}$.

A state $s$ of $D$ is a labelling of each crossing of $D$ by
either an A-chord or a B-chord. We write $a(s)$ (resp. $b(s)$) for
the number of A-chords (resp. B-chords) of the state $s$, and
$|s|$ for the number of components of the diagram $sD$, which is
$D$ after the $s$-smoothing of $D$. Precisely $sD$ is obtained
smoothing every crossing in $D$ according to the type of chord
associated to the crossing by the state, as shown in Figure 1. We
will draw a small chord with the letters $A$ or $B$ to remember
which was the state. In this way we can reconstruct the diagram
$D$ from $sD$ and the chords, by just reversing the smoothings
shown in Figure 1. An example is shown in Figure 2.
\begin{figure}[h]
  \begin{center}
    \includegraphics[height=2cm, width=12cm]{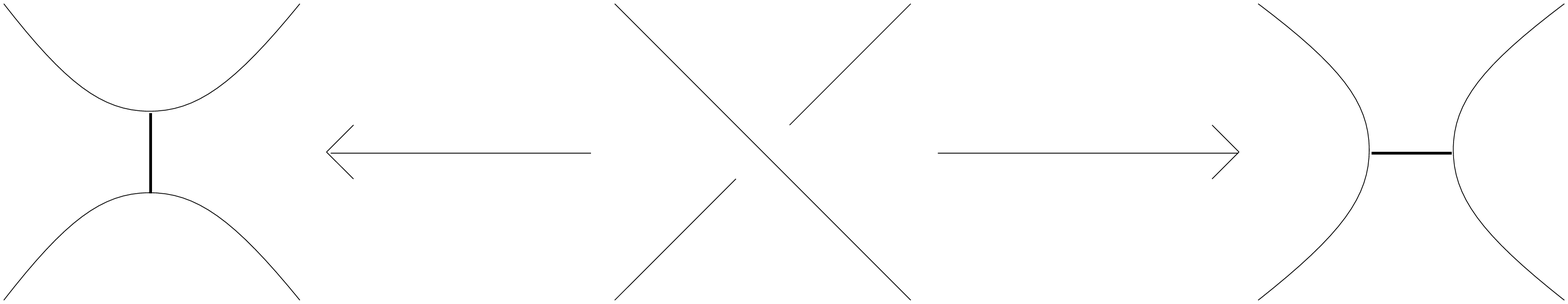}
    \caption{$s_A$ and $s_B$-smoothing of a crossing and corresponding chords.}
      \begin{picture}(0,0)
        \put(-135,61){$A$}
        \put(133,68){$B$}
        \put(37,72){{\small $s_B$-smoothing}}
        \put(-94,72){{\small $s_A$-smoothing}}
      \end{picture}
  \end{center}
\end{figure}
\begin{figure}[h]
  \begin{center}
    \includegraphics[height=2cm, width=12cm]{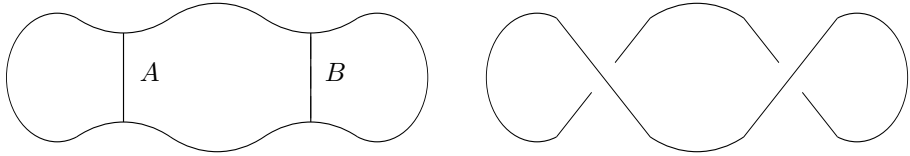}
    \caption{From $sD$ and the chords to $D$.}
      \begin{picture}(0,0)
        \put(-50,61){$B$}
        \put(-120,61){$A$}
      \end{picture}
  \end{center}
\end{figure}
For a state $s$ we denote by $max(s)$ (resp. $min(s)$) the highest
(resp. the lowest) degree of $\langle D,s\rangle$. The extreme
states $s_A$ and $s_B$ are defined by the equalities $a(s_A)=c(D)$
and $b(s_B)=c(D)$ respectively, where $c(D)$ is the number of
crossings of the diagram $D$. Write $m=min(s_B)$ and $M=max(s_A)$.
Clearly $m=-c(D)-2|s_B|+2$ and $M=c(D)+2|s_A|-2$. These numbers
$m$ and $M$ will be called the extreme states degrees of $\langle
D \rangle $ and their corresponding coefficients $a_m$ and $a_M$
in $\langle D \rangle$ will be called the extreme states
coefficients of $\langle D \rangle $. It turns out that (see
Proposition 1)

$$\langle D \rangle =a_mA^m+a_{m+4}A^{m+4}+\dots
+a_{M-4}A^{M-4}+a_MA^{M}.$$

In this paper we deal with the question of finding arbitrary
extreme coefficients for $\langle D \rangle $. Two different
approachings will be given:

\vspace{0.3cm}

The first approaching follows Yonju Bae and Morton \cite{main}. In
their paper a connection between $a_M$ and graph theory is given.
The coefficient $a_M$ appears to be, up to sign, the value
$f(G_D)$ of a certain graph $G_D$, constructed from $s_AD$ and
called the Lando's graph of $sD$. In general, for any graph $G$,
$f(G)=\sum_C(-1)^{|C|}$ where $C$ runs over all the independent
subsets of vertices of $G$, where independent means that there is
no edge joining two vertices of $C$. In \cite{main} it arises the
question of if any integer can be realized as $f(G)$ for some
graph $G_D$. We have found a positive answer to this question. In
fact, we will exhibit examples (in fact complete families) of
graphs, and from these we will reconstruct prime knot diagrams
with arbitrary extreme states coefficients $a_m$ and $a_M$. This
will be done in the second section. The third section will be
dedicated to explain further constructions of graphs, which in
many cases give easier examples of knots with the wanted extreme
states coefficients. Here ``easier'' means much fewer crossings.

\vspace{0.3cm}

The second approaching is treated in the fourth and final section.
Here we exhibit examples of prime knot diagrams for which the
extreme states coefficients $a_m$ and $a_M$ are zero, and the next
coefficients $a_{m+4}$ and $a_{M-4}$ take arbitrary values. The
idea is just to look at the more general circle graph defined by
$sD$, rather than the Lando's graph $G_D$, and use a very simple
trick for counting $a_{M-4}$ in special cases. We do not know if
there is a complete nice interpretation of $a_{M-4}$ in terms of
graph theory, parallel to that one in which $a_M$ is described in
terms of $G_D$ in the second section. A really more interesting
question is the following: is it possible to get any extreme
coefficient when the spread of the Jones polynomial is previously
fixed?

\vspace{0.3cm}

\section{Extreme states coefficients of Jones polynomials and graph
theory.}

We begin by recalling some very basic facts about the Kauffman
bracket $\langle D \rangle$ of an unoriented diagram $D$.

\begin{proposition} \label{well-known}

\noindent \ \ (i) All degrees in $\langle D \rangle$ are congruent
modulo four.

\vspace{0.3cm}

\noindent \ (ii) $max(s)\leq M=max(s_A)$ with equality if and only
if $|s|=|s_A|+b(s)$.

\vspace{0.3cm}

\noindent (iii)  $min(s)\geq m=min(s_B)$ with equality if and only
if $|s|=|s_B|+a(s)$.

\vspace{0.3cm}

\noindent (iv) The highest (resp. lowest) degree of $\langle D
\rangle$ is less (resp. great) or equal than $M$ (resp. $m$).

\vspace{0.3cm}

\noindent (v) A state $s$ contributes to $a_M$ if and only if
$s\in\Gamma _A=\{ s/|s|=|s_A|+b(s) \}$. The contribution of
$s\in\Gamma _A$ to $a_M$ is $(-1)^{|s_A|-1}(-1)^{b(s)}$.

\vspace{0.3cm}

\noindent (vi) A state $s$ contributes to $a_m$ if and only if
$s\in\Gamma _B=\{ s/|s|=|s_B|+a(s) \}$. The contribution of
$s\in\Gamma _B$ to $a_m$ is $(-1)^{|s_B|-1}(-1)^{a(s)}$.

\vspace{0.3cm}

\noindent (vii) $a_M=(-1)^{|s_A|-1}\sum _{s\in \Gamma
_A}(-1)^{b(s)}$ and $a_m=(-1)^{|s_B|-1}\sum _{s\in \Gamma
_B}(-1)^{a(s)}$.
\end{proposition}
\begin{proof}
We calculate the difference
\begin{eqnarray*}
max(s_A)-max(s) &=&c(D)+2|s_A|-2-a(s)+b(s)-2|s|+2\\
&=&2b(s)+2|s_A|-2|s|
\end{eqnarray*}
Now the key point is that if two states $s$ and $s'$ differ in the
label of only one crossing, then either $|s'|=|s|+1$ or
$|s'|=|s|-1$, depending on whenever the two strings that appear
after $s$-smoothing the crossing belong or not to the same
component of $sD$. It follows that for an arbitrary state $s$ we
have that there are non-negative integers $p$ and $n$ such that
$b(s)=p+n$ and $|s|=|s_A|+p-n$. Then
$b(s)+|s_A|-|s|=p+n+|s_A|-|s_A|-p+n=2n$ is an even number great or
equal than zero. The completation of the proof is left to the
reader.
\end{proof}

\vspace{0.3cm}

\noindent {\bf Remark} \label{remark} In order to find $sD$, one
can start with $s_AD$, choose an order in the set of the $b(s)$
crossings labelled with a B-chord in the state $s$ and perform the
opposite smoothing in these crossings following this order. In
this way, one associate $+1$ (resp. $-1$) to each one of the
$b(s)$ crossings if after the $s_B$-smoothing of this crossing we
get one more (resp. fewer) component. By definition $p$ (resp.
$n$) is the number of associated $+1$ (resp. $-1$). This
non-negative integers does not depend on the chosen order of the
$b(s)$ crossings, although the association of $+1$ and $-1$ to
each one of the $b(s)$ crossings does.

\vspace{0.3cm}

Now we recall the connection between $a_M$ and graph theory given
in \cite{main}. In order to obtain the Lando's graph of $D$ start
with $s_AD$ and delete the A-chords joining two different
components. In this way we obtain a bipartite circle graph (BCG).
Now define the Lando's graph of $D$ as the graph that appears
taking a vertex for every A-chord of the BCG, and joining two
vertices with an edge if and only if the endpoints of the
corresponding A-chords in the BCG alternate in the same component
of the BCG. An example of $s_AD$, the corresponding BCG and the
Lando's graph is shown in Figure 3.
\begin{figure}[h]
  \begin{center}
    \includegraphics[height=3cm, width=12cm]{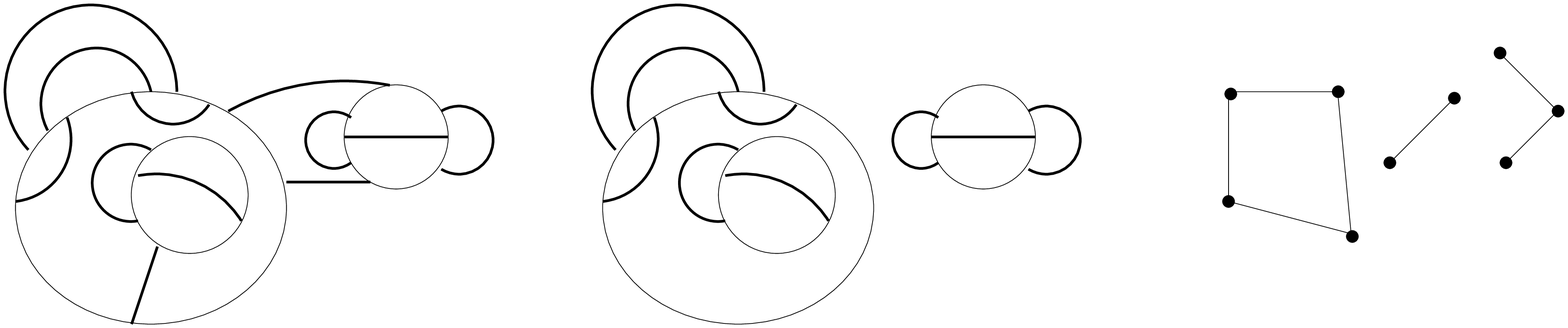}
    \caption{$s_AD$, the BCG and the Lando's graph $G_D$.}
      \begin{picture}(0,0)
        \put(-110,37){$s_AD$}
        \put(17,37){BCG}
        \put(114,37){$G_D$}
      \end{picture}
  \end{center}
\end{figure}
A subset $C$ of vertices of a graph $G$ is said to be independent
if there is no edge in $G$ joining two vertices of $C$. We define
$f(G)=\sum_C(-1)^{|C|}$ where $C$ runs over the independent sets
of vertices of $G$.

Then
\begin{eqnarray*}
a_M &=&(-1)^{|s_A|-1}\sum _{s\in\Gamma _A}(-1)^{b(s)}\\
&=&(-1)^{|s_A|-1}\sum _{C}(-1)^{|C|}\\ &=&(-1)^{|s_A|-1}f(G_D).
\end{eqnarray*}
The first equality is given by Proposition~\ref{well-known} and
the third equality is just the definition of $f(G_D)$. In order to
check the second equality, we think of a state $s$ as the set of
$b(s)$ A-chords of $s_AD$ which correspond to the $b(s)$ B-chords
of $s$. Then we have that $s\in \Gamma _A$ if and only if the two
following conditions occur:

(1) The endpoints of every A-chord lies in the same component.

(2) The endpoints of two A-chords lying in the same component do
not alternate.

The example in Figure 4 is exhibited in \cite{main}, where the
corresponding Lando's graph $G_D$ has $f(G_D)=3$.
\begin{figure}[h]
  \begin{center}
    \includegraphics[height=3cm, width=12cm]{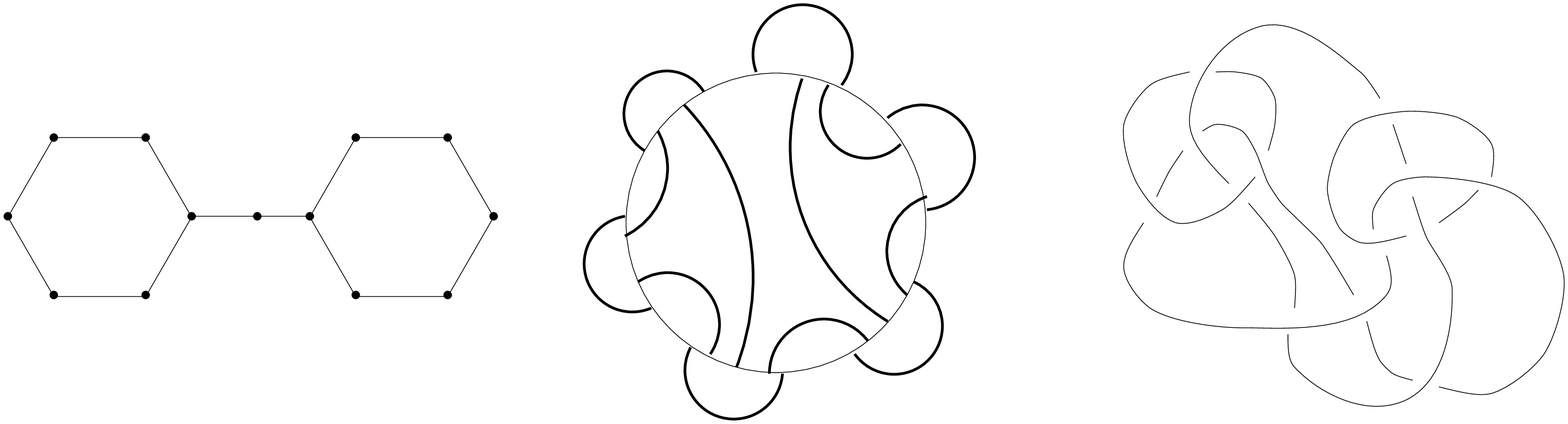}
    \caption{A Lando's graph $G_D$ with $f(G_D)=3$ and the corresponding BCG and link diagram.}
      \begin{picture}(0,0)
        \put(-120,66){$G_D$}
        \put(0,40){BCG}
        \put(105,42){$D$}
      \end{picture}
  \end{center}
\end{figure}
The question formulated in \cite{main} is if any integer $n$ can
be realized as $f(G_D)$ for a graph $G_D$ arising from a diagram
$D$. Since these graphs arise from BCG, we will call them ``graphs
convertible in BCG''.

\vspace{0.3cm}

\noindent {\bf Some graph theory.}

\vspace{0.3cm}

Since the graphs $G_D$ must be always convertible in BCG, the
examples shown in Figure 5 are not allowed. On the other hand the
graphs in consideration are not necessarily planar graphs. Figure
6 exhibits such an example.
\begin{figure}[h]
  \begin{center}
    \includegraphics[height=1cm, width=6cm]{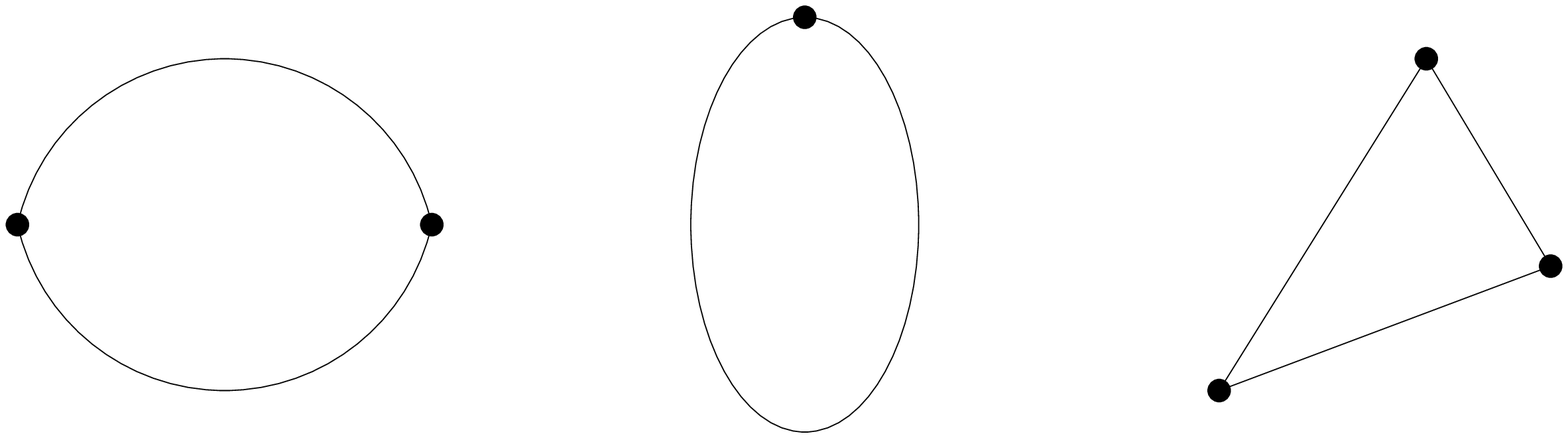}
    \caption{Examples of graphs non-convertible in BCG.}
  \end{center}
\end{figure}
\begin{figure}[h]
  \begin{center}
    \includegraphics[height=2cm, width=5cm]{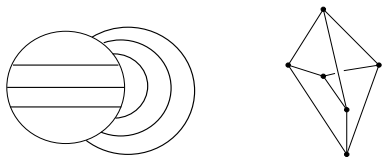}
    \caption{Non-planar Lando's graph.}
  \end{center}
\end{figure}

Calculation of $f(G)$ can be simplified using some readily
established properties, as explained in \cite{main}. Here it is a
very brief description of these properties:

\vspace{0.3cm}

\noindent $\bullet$ Law of recursion. If $G$ is a graph and $v$ is
a vertex of $G$, we will denote by $G-v$ the graph obtained from
$G$ by deleting the vertex $v$ and its incident edges. Let $\{
v_1, \ldots ,v_k \}$ be the set of the neighbour vertices of $v$
in $G$ (by definition these are the vertices of $G$ joined to $v$
by an edge). We will denote by $G-Nv$ the graph $(\ldots
((G-v)-v_1)\ldots )-v_k$. Then the law of recursion says that
$f(G)=f(G-v)-f(G-Nv)$.

\vspace{0.3cm}

\noindent $\bullet$ Law of multiplication. If a graph $G$ is the
disjoint union of two graphs $G_1$ and $G_2$, then
$f(G)=f(G_1)(G_2)$.

\vspace{0.3cm}

\noindent $\bullet$ Law of duplication. Suppose that $v$ and $w$
are two non-neighbour vertices of a graph $G$, and the set of
neighbour vertices of $v$ is included in the set of neighbour
vertices of $w$. Then we say that $G$ is a duplication of $G-w$
and we have $f(G)=f(G-w)$.

\vspace{0.3cm}

We consider examples that will be used later. Let $L_n$ be the
graph with $n$ vertices shown in Figure 7. Clearly $f(L_2)=-1$. By
using duplication in the third vertex and then multiplication we
get the formula $f(L_n)=-f(L_{n-3})$.
\begin{figure}[h]
  \begin{center}
    \includegraphics[height=0.2cm, width=8cm]{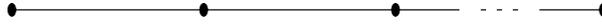}
    \caption{Graph $L_n$.}
  \end{center}
\end{figure}
Let $C_n$ be the polygon with $n$ vertices, $n\geq 3$. By
recursion we obtain $f(C_n)=f(L_{n-1})-f(L_{n-3})$. In particular,
for the hexagon $H=C_6$ we have $f(H)=2$. We now explain how to
find different families $\{G_r\} _{r\in {\bf Z}}$ of graphs with
$f(G_r)=r$ for any integer $r$. From now on, we will write $G^v$
to mean the pair $(G,v)$, where $G$ is a graph and $v$ is a
particular vertex of $G$.

\vspace{0.3cm}

\noindent {\bf Definition} Let $G$ be a planar graph convertible
in BCG and let $v$ be a vertex of $G$. We say that $G^v$ if a
brick of type $(n,k)$ if $f(G)=n$ and $f(G-v)=k$.

\vspace{0.3cm}

Our main example will be the hexagon $H$ with any arbitrary vertex
chosen, which is a brick of type $(2,1)$. It is a graph
convertible in BCG, as we show in Figure 8.
\begin{figure}[h]
  \begin{center}
    \includegraphics[height=2cm, width=5cm]{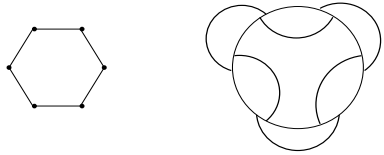}
    \caption{The hexagon $H$ is convertible in BCG.}
  \end{center}
\end{figure}
Now we describe the basic construction that will be used for
giving our examples: Let $G_1$ and $G_2$ be two graphs and $v_1$
and $v_2$ two vertices of $G_1$ and $G_2$ respectively. Then
$G_1^{v_1}*G_2^{v_2}$ will denote the new graph obtained from the
disjoint union of $G_1$ and $G_2$ by joining $v_1$ and $v_2$ with
an extra edge (see Figure 9). Note that $G_1^{v_1}*G_2^{v_2}$ is
convertible in BCG if both $G_1$ and $G_2$ are.
\begin{figure}[h]
  \begin{center}
    \includegraphics[height=2cm, width=7cm]{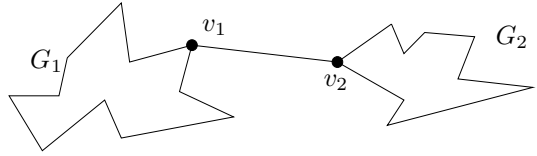}
    \caption{Graph $G_1^{v_1}*G_2^{v_2}$.}
       \begin{picture}(0,0)
         \put(-91,66){$G_1$}
         \put(-26,80){$v_1$}
         \put(20,59){$v_2$}
         \put(85,75){$G_2$}
       \end{picture}
  \end{center}
\end{figure}
\begin{lemma} \label{lemmaprime}
Let $G$ be a graph and $v$ a vertex of $G$. If $G^v$ is a brick of
type $(n,k)$, then $(G',v')$ is a brick of type $(n+k,k)$, where
$G'=G^v*H^w$ and $v'$ is a vertex of $H$ adjacent to $w$ (see
Figure 10).
\end{lemma}
\begin{figure}[h]
  \begin{center}
    \includegraphics[height=2cm, width=6cm]{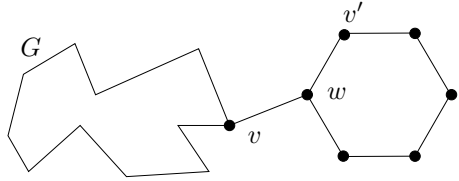}
    \caption{Graph $(G',v')$ arising from $(G,v)$.}
       \begin{picture}(0,0)
         \put(-80,80){$G$}
         \put(6,47){$v$}
         \put(36,64){$w$}
         \put(42,93){$v'$}
       \end{picture}
  \end{center}
\end{figure}
\noindent
\begin{proof}
We have to show that $f(G')=n+k$ and $f(G'-v')=k$.

\vspace{0.3cm}

We have $f(G'-v')=(f(L_2))^2f(G-v)=(-1)^2k=k$ where we have used
duplication, first in the vertex $v_1$ and then in the vertex $v$,
multiplication and the fact that $f(L_2)=-1$ (see Figure 11).
\begin{figure}[h]
  \begin{center}
    \includegraphics[height=2cm, width=6cm]{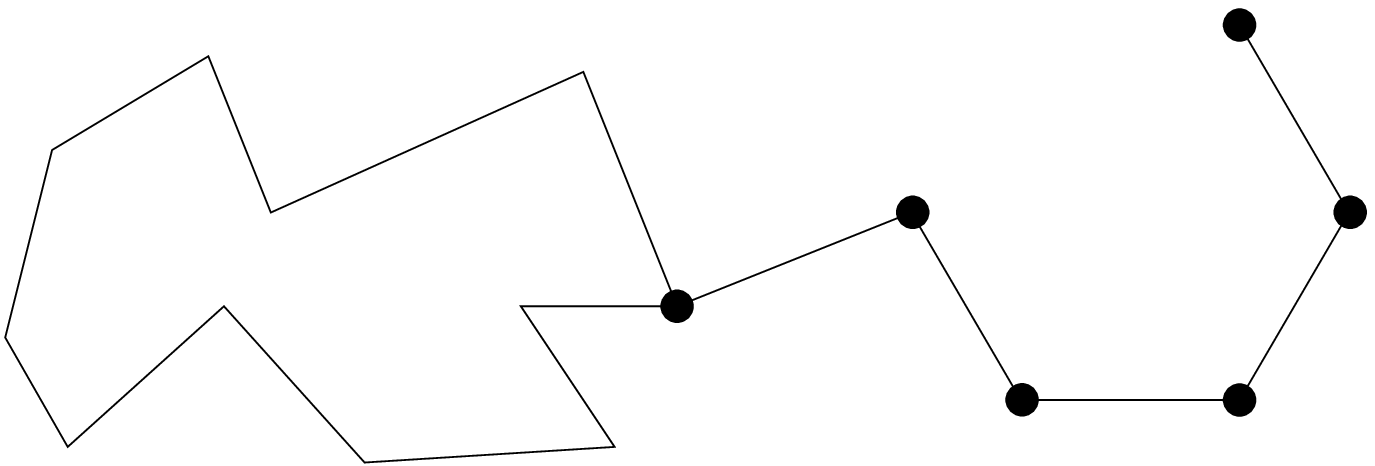}
    \caption{$G'-v'$.}
       \begin{picture}(0,0)
         \put(5,48){$v$}
         \put(75,40){$v_1$}
       \end{picture}
  \end{center}
\end{figure}

On the other hand
\begin{eqnarray*}
f(G') &=&f(G'-v')-f(G'-Nv') \ \ \ \ (recursion)\\
 &=&k-f(L_3)f(G') \ \ \ \ \ \ \ \ \ \ \ \ \ \ \ \ \ (multiplication)\\
 &=&k-(-1)n \ \ \ \ \ \ \ \ \ \ \ \ \ \ \ \ \ \ \ \ \ \ \ \ (f(L_3)=-1)\\
 &=&n+k.
\end{eqnarray*}

\vspace{0.3cm}

Finally, $G'$ is convertible in BCG since $G$ and $H$ are.
\end{proof}

\vspace{0.3cm}

\begin{lemma} \label{lemmaminus}
Let $G$ be a graph. Then $f(G^-)=-f(G)$ where $G^-=G^v*L_3^w$, $v$
being any vertex of $G$ and $w$ being an extreme vertex of $L_3$.
Moreover, $G^-$ is convertible in BCG if $G$ is (see Figure 12).
\end{lemma}
\begin{figure}[h]
  \begin{center}
    \includegraphics[height=2cm, width=6cm]{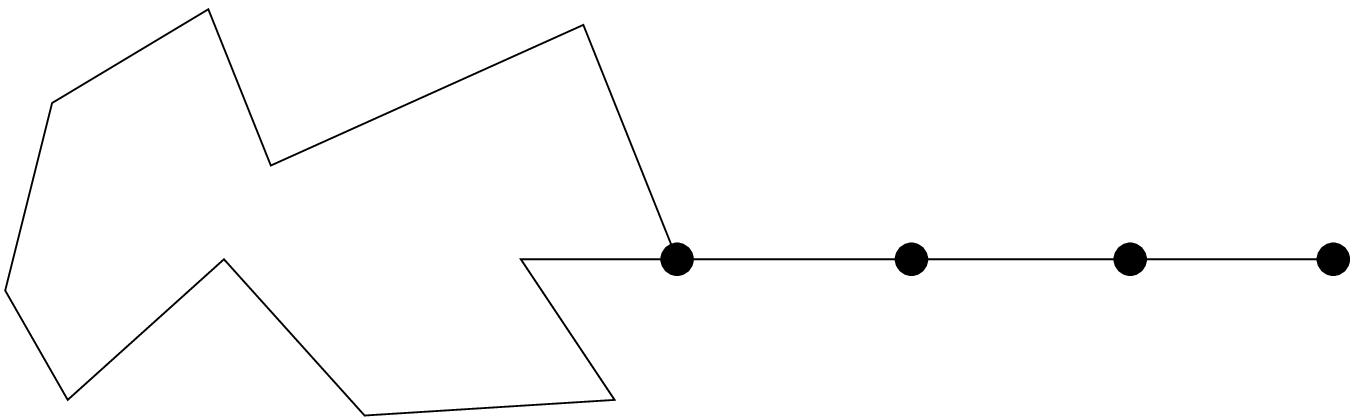}
    \caption{Graph $G^-$.}
       \begin{picture}(0,0)
         \put(-84,81){$G$}
         \put(5,49){$v$}
       \end{picture}
  \end{center}
\end{figure}
\noindent \begin{proof} We have
\begin{eqnarray*}
f(G^-)&=&f(G-w) \ \ \ \ \ \ \ \ \ \ \ \ (duplication)\\
&=&f(L_2)f(G) \ \ \ \ \ \ \ \ \ \ (multiplication)\\ &=&-f(G) \ \
\ \ \ \ \ \ \ \ \ \ \ \ \ \ (f(L_2)=-1).
\end{eqnarray*}

\vspace{0.3cm}

\noindent Finally, $G^-$ is convertible in BCG since $G$ and $L_3$
are.
\end{proof}

\vspace{0.3cm}

\begin{theorem} \label{main} For any integer $n$
there is a planar graph $G_{n-1}$ convertible in bipartite circle
graph such that $f(G_{n-1})=n$.
\end{theorem}
\noindent
\begin{proof}Consider the brick $H^v$ of type $(2,1)$ and apply
r times Lemma 1 to get a planar graph $G_{r+1}$ convertible in BCG
such that $f(G_{r+1})=r+2$ (see Figure 13).

\vspace{0.3cm}

On the other hand $L_4$, $L_1$ and $L_2$ are planar graphs
convertible in BCG with $f(L_4)=1$, $f(L_1)=0$ and $f(L_2)=-1$
respectively.

\vspace{0.3cm}

Finally, for all integer $r\geq 1$ we have that $G_{r+1}^-$ is a
planar graph convertible in BCG with $f(G_{r+1}^-)=-(r+2)$
according to Lemma~\ref{lemmaminus} (see Figure 14).
\end{proof}
\begin{figure}[h]
  \begin{center}
    \includegraphics[height=2cm, width=10cm]{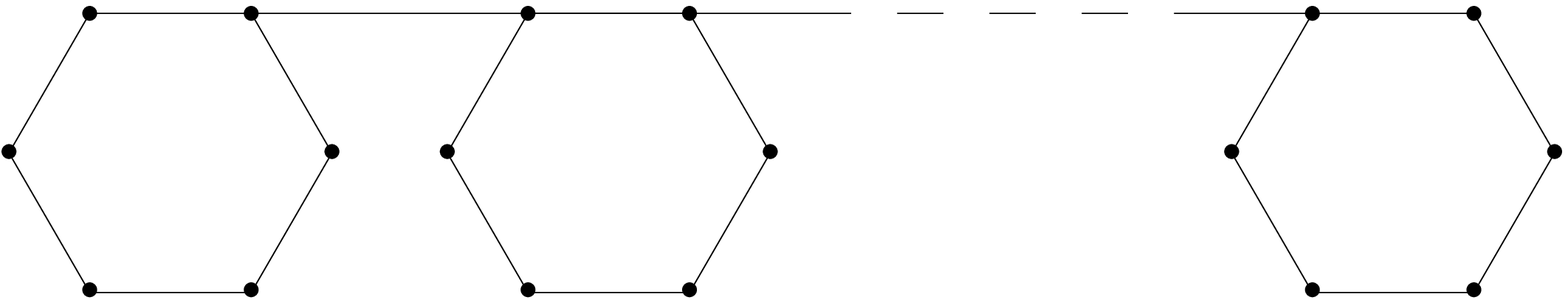}
    \caption{The graph $G_r$.}
       \begin{picture}(0,0)
         \put(20,79){{\small $r$ hexagons}}
       \end{picture}
  \end{center}
\end{figure}
\begin{figure}[h]
  \begin{center}
    \includegraphics[height=2cm, width=12cm]{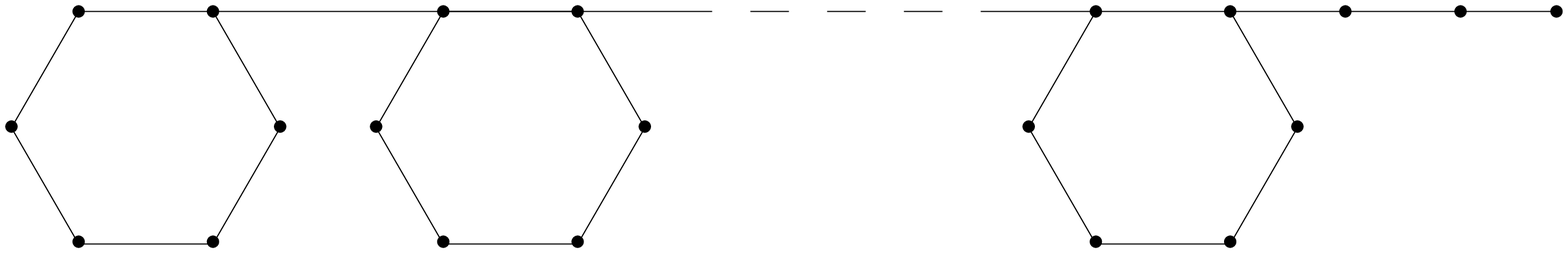}
    \caption{The graph $G_r^-$.}
        \begin{picture}(0,0)
             \put(-7,79){{\small $r$ hexagons}}
        \end{picture}
  \end{center}
\end{figure}

\vspace{0.3cm}

We will now construct a prime knot diagram with arbitrary extreme
states coefficients, starting with the graphs $G_r$. This process
is illustrated in figures 15, 16, 17 and 18. We first reconstruct
the associated bipartite circle graph (Figure 15). Changing every
A-chord $\raisebox{-.8mm}{\epsfysize.15in \epsffile{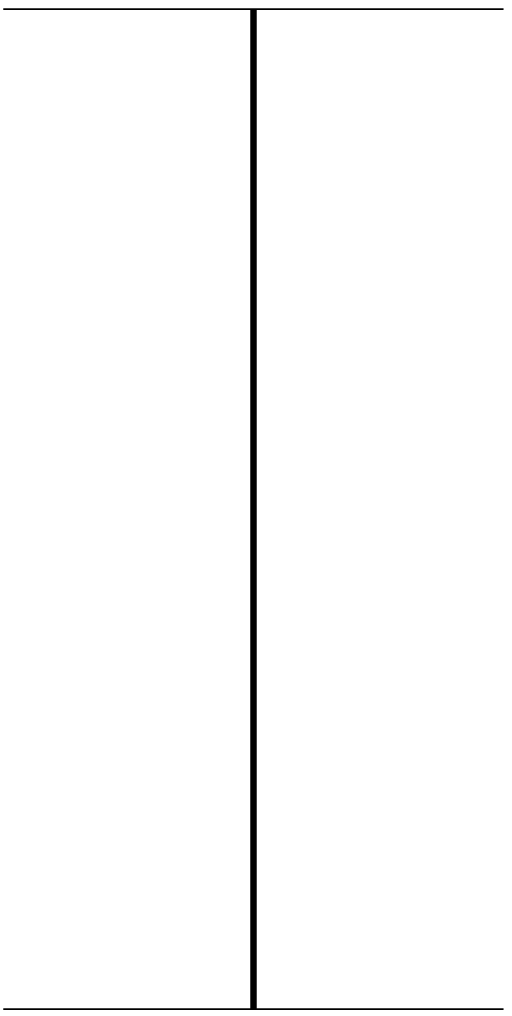}}$
by a crossing $\raisebox{-.8mm}{\epsfysize.15in
\epsffile{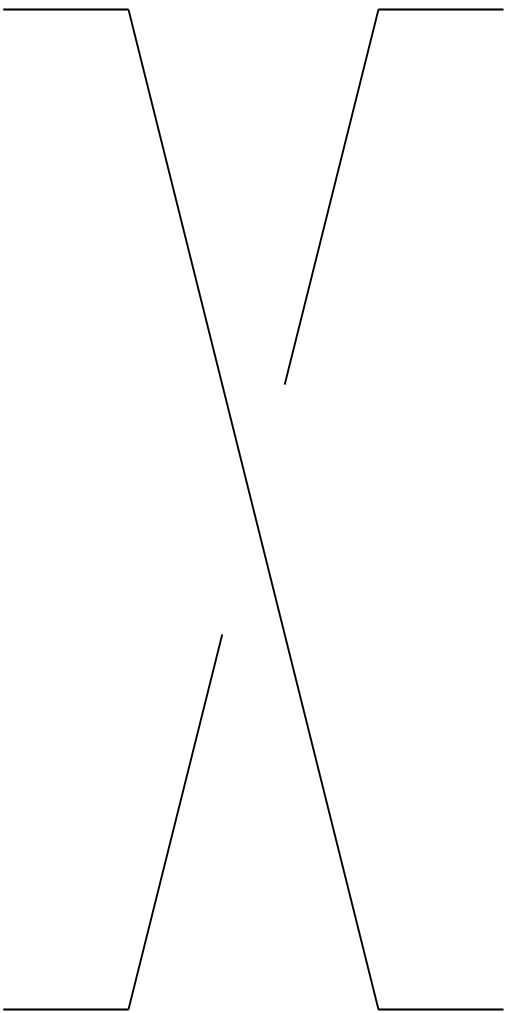}}$ we see that this BCG is $s_AD_r'$ where
the ($3-$components) diagram $D_r'$ is shown in Figure 16. Figure
17 shows $s_BD_r'$, proving that $D_r'$ is quite far from being a
-adequate diagram (the coefficient $a_m$ of $\langle D_r'\rangle $
is not $\pm 1$ in general since there are other states in $\Gamma
_B$ apart from $s_B$). Because of this we modify $D_r'$ to produce
the diagram $D_r$ in Figure 18.
\begin{figure}[h]
  \begin{center}
    \includegraphics[height=3cm, width=12cm]{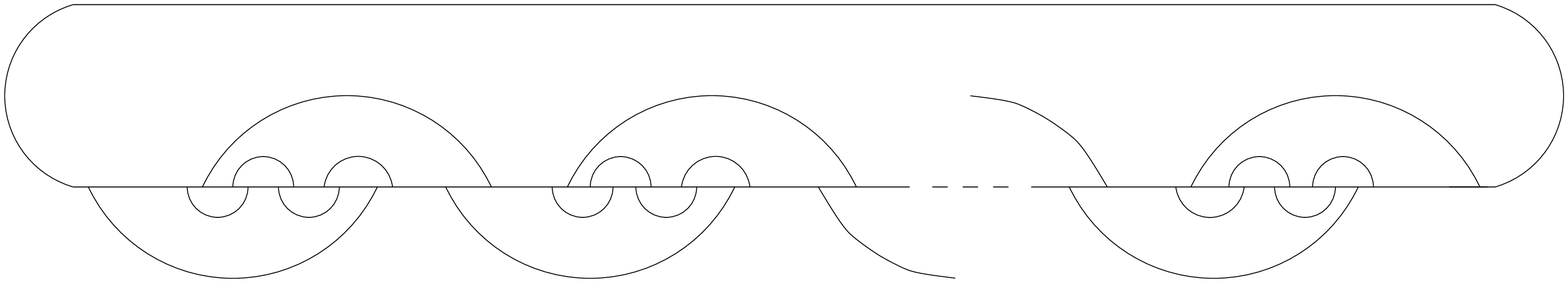}
    \caption{The BCG obtained from $G_r$.}
  \end{center}
\end{figure}
\begin{figure}[h]
  \begin{center}
    \includegraphics[height=3cm, width=12cm]{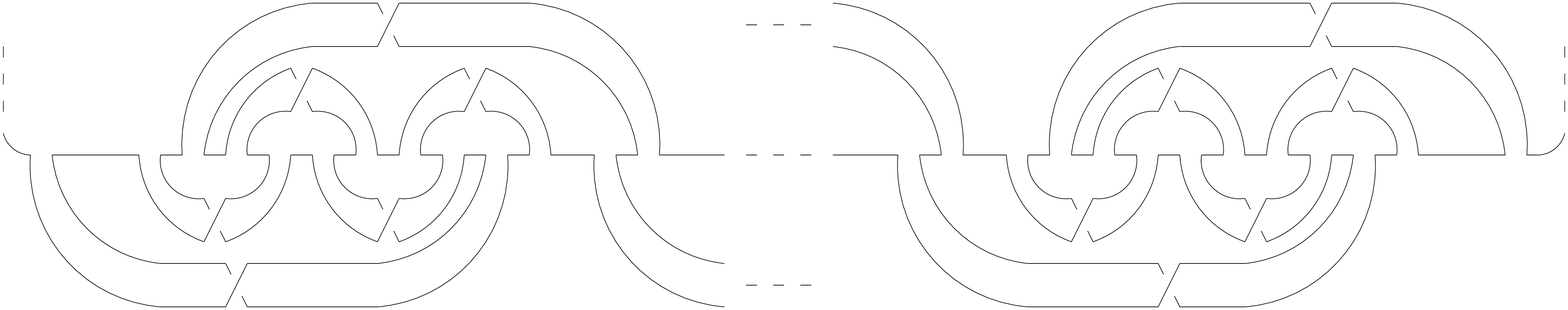}
    \caption{The diagram $D_r'$.}
  \end{center}
\end{figure}
\begin{figure}[h]
  \begin{center}
    \includegraphics[height=3cm, width=12cm]{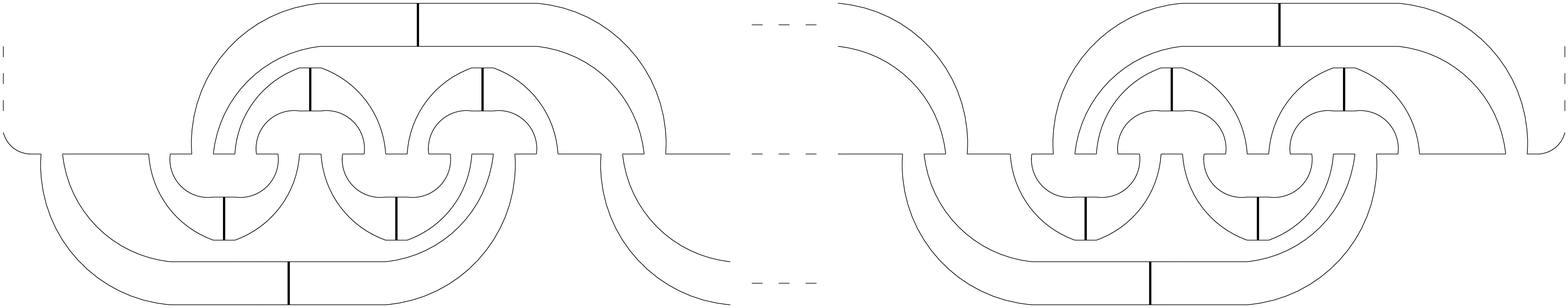}
    \caption{$s_BD_r'$: $D_r'$ is not -adequate.}
  \end{center}
\end{figure}
\begin{figure}[h]
  \begin{center}
    \includegraphics[height=3cm, width=12cm]{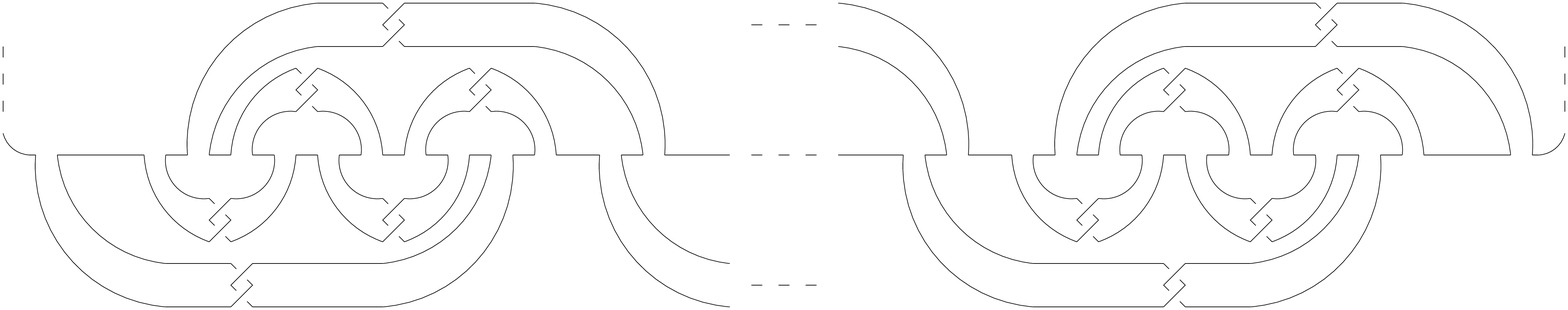}
    \caption{The diagram $D_r$.}
  \end{center}
\end{figure}

\begin{theorem}\label{Dr} $\langle D_r \rangle
=A^m+\dots + (r+1)A^M$ where $m=-24r-4$ and $M=12r$. In particular
$span(\langle D_r \rangle )=4(9r+1)$. Moreover, $D_r$ is a prime
knot diagram with $12r$ crossings.
\end{theorem}
\noindent
\begin{proof}First we fix our attention in $D_r$ after
$s_A$-smoothing. Note that $|s_AD_r|=1$ hence $M=c(D)=12r$. On the
other hand the Lando's graph of $D_r$ appear to be a duplication
of $G_r$, hence $a_M=(-1)^{|s_A|-1}f(G_r)=r+1$.

\vspace{0.3cm}

Now apply the $s_B$-smoothing to $D$. We have $|s_BD_r|=6r+3$
($6r$ components are given by the small circles
$\raisebox{-.8mm}{\epsfysize.15in \epsffile{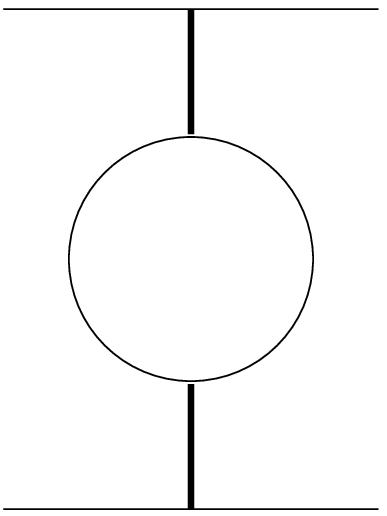}}$ and
the other three components are those appearing in Figure 17),
hence $m=-12r-2(6r+3)+2=-24r-4$. On the other hand the Lando's
graph (respect to the $s_B$-smoothing) is the empty set, hence
$a_m=(-1)^{|s_B|-1}f(\emptyset )=1$.
\end{proof}

\vspace{0.3cm}

We now join in a very particular way $D_r$ and $\bar D_s$ (the
mirror image of $D_s$) in order to control simultaneously both
extreme states coefficient (we avoid the obvious solution given by
the appropriate connected sum in order to get a prime diagram):

\begin{theorem}\label{Drs} Let $D_{rs}$ be the diagram
shown in Figure 19. Then $\langle D_{rs} \rangle =(s+1)A^m+\dots
+(r+1)A^M$ where $m=-24r-12s-6$ and $M=12r+24s+6$. In particular,
$span(\langle D_{rs} \rangle )=36(r+s)+12$. Moreover, $D_{rs}$ is
a prime knot diagram with $12(r+s)+2$ crossings.
\end{theorem}
\begin{figure}[h]
  \begin{center}
    \includegraphics[height=5cm, width=12cm]{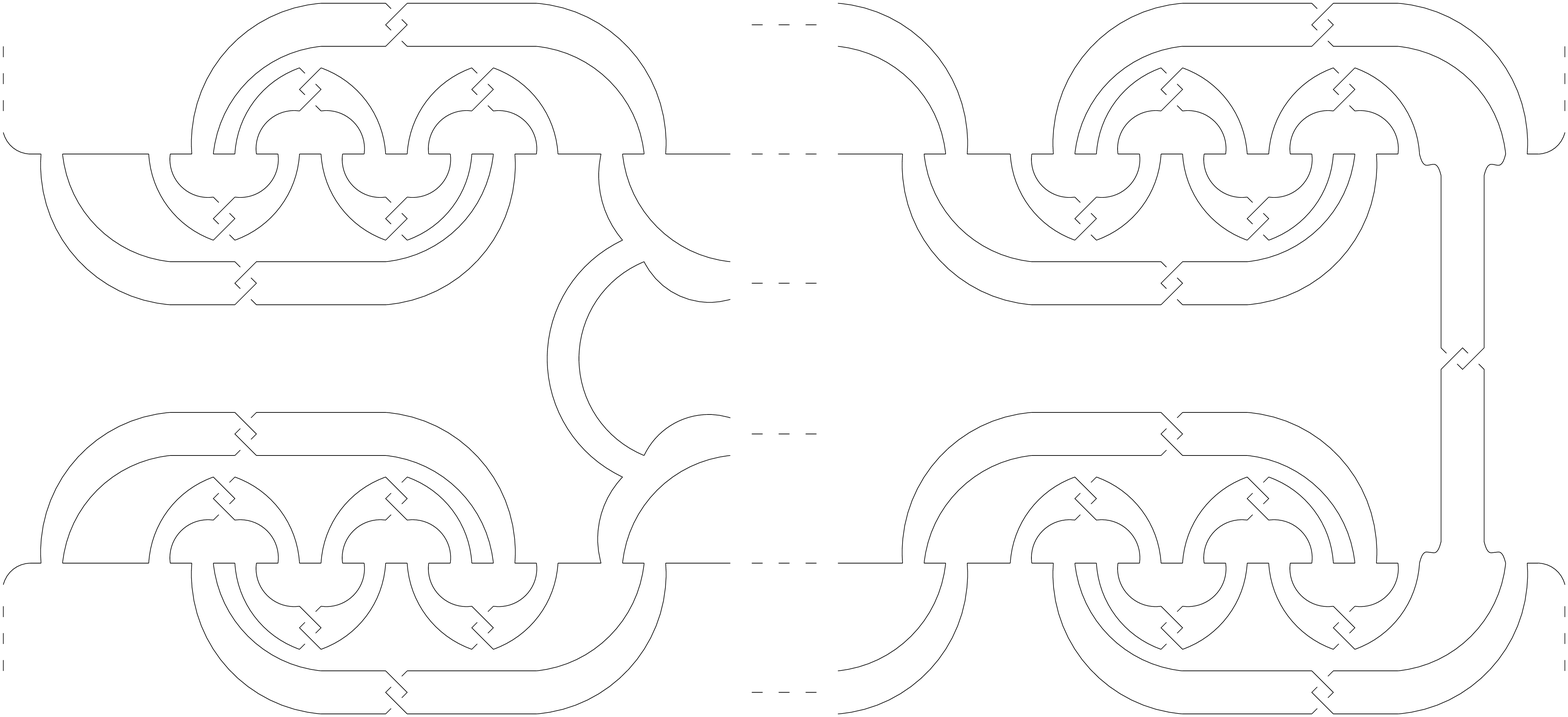}
    \caption{The diagram $D_{rs}$.}
  \end{center}
\end{figure}
\noindent
\begin{proof}Note first that $c(D_{rs})=c(D_r)+c(\bar D_s)+2=12r+12s+2$.

\vspace{0.3cm}

Now, the very special way in which we join $D_r$ and the mirror
image $\bar D_s$ of $D_s$ gives the equalities $|s_A|=|s_A(\bar
D_s)|=|s_B(D_s)|=6s+3$ and the fact that the Lando's graph is
still a duplication of $G_r$. It follows that $M=12r+24s+6$ and
$a_M=(-1)^{|s_A|-1}f(G_r)=r+1$.

\vspace{0.3cm}

\noindent Analogously, $|s_B|=|s_B(D_r)|=6r+3$ and the Lando's
graph is still a duplication of $G_s$, hence $m=-24r-12s-6$ and
$a_m=(-1)^{|s_B|-1}f(G_s)=s+1$.~\end{proof}

\vspace{0.3cm}

A small refinement of the last result allows us to modify the
signs of the extreme states coefficients:
\begin{theorem}\label{Drsalpha} Let $\alpha $
be an odd integer great than $1$. Let $D_{rs}^\alpha$ be the
diagram $D_{rs}$ with a modification in the way in which $D_r$ and
$\bar D_s$ are joined on the left, as shown in Figure 20. Then
$\langle D_{rs}^\alpha \rangle =-(s+1)A^m+\dots +(r+1)A^M$ where
$m=-24r-12s-\alpha -8$ and $M=12r+24s+3\alpha +4$. In particular
$span(\langle D_{rs}^\alpha \rangle )=36(r+s)+4\alpha +12$.
Moreover, $D_{rs}^\alpha$ is a prime knot diagram with
$12(r+s)+2+\alpha $ crossings.
\begin{figure}[h]
  \begin{center}
    \includegraphics[height=1cm, width=4cm]{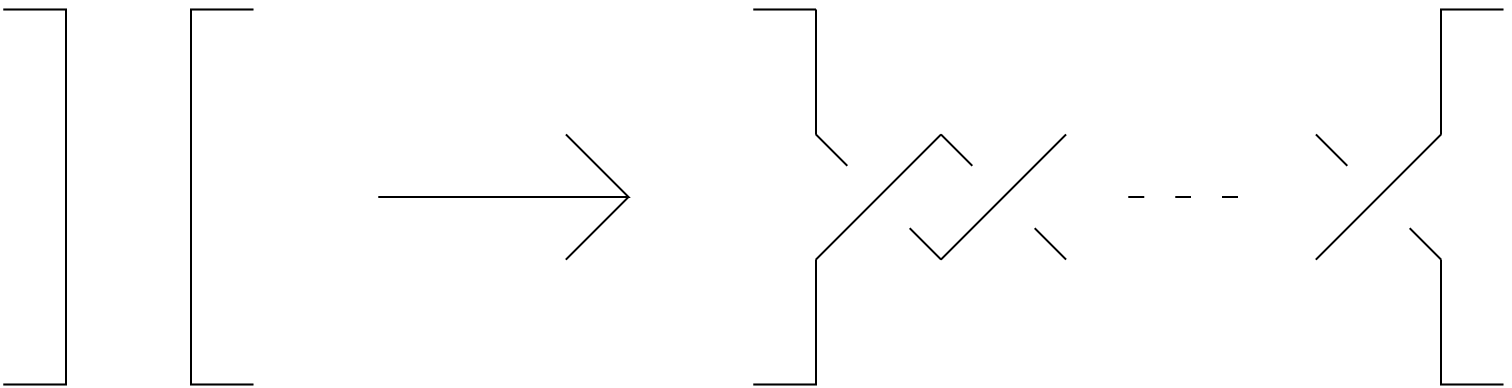}
    \caption{A small modification of $D_{rs}$ produces $D_{rs}^\alpha $.}
        \begin{picture}(0,0)
             \put(30,38){$\alpha $}
        \end{picture}
  \end{center}
\end{figure}
\end{theorem}
\noindent
\begin{proof} We have that $c(D_{rs}^\alpha)=c(D_{rs})+\alpha $,
$|s_A|=|s_A(D_{rs})|+\alpha -1=6s+\alpha +2$ and
$|s_B|=|s_B(D_{rs})|+1=6r+4$. In addition the $s_A$ and $s_B$
Lando's graphs are still the duplications of $G_r$ and $G_s$
respectively.
\end{proof}

\vspace{0.3cm}

\section{More graph theory.}

Let us consider the specific example provided by Theorem \ref{Dr}
for the extreme state coefficient $a_M=41$. The knot diagram
$D_{40}$ has $12 \times 40=480$ crossings! If we take the non
-adequate link diagram with three components $D_{40}'$, we still
have $6\times 40=240$ crossings! In this section we develop other
graph constructions, providing easier examples for many possible
extreme states coefficients. Here ``easier'' means always that the
diagram has fewer crossings. In many cases, easier means as well
that the corresponding span is lower.

\vspace{0.3cm}

\noindent {\bf Building with bricks}

\vspace{0.3cm}

Let $\{ G_1^{v_1},\dots ,G_k^{v_k} \}$ a set of graphs with an a
chosen vertex $v_i\in G_i$ for every $1\leq i \leq k$. We denote
by $S=S(G_1^{v_1},\dots ,G_k^{v_k})$ the graph shown in Figure
21a, and we call it ``simple building'' constructed with the
bricks $G_1^{v_1},\dots ,G_k^{v_k}$. Precisely, $S$ can be defined
in steps using the operation $*$ introduced in the second section.
Let $W$ be the only vertex of $L_1$. Then

$$S= ((\dots ((L_1^w*G_1^{v_1})^w*G_2^{v_2})^w*\dots
)^w)*G_k^{v_k}.$$

Other related construction can be obtained bay introducing $k$
extra vertices $w_1, \dots ,w_k$, one in every edge joining $v_i$
and $w$ in the graph $S$. The new graph is denoted by
$C=C(G_1^{v_1},\dots ,G_k^{v_k})$ and it is shown in Figure 21b.
We will call it ``complicated building'' constructed with the
bricks $G_1^{v_1},\dots ,G_k^{v_k}$. Its precise definition using
the operation $*$ is left to the reader. In both constructions the
vertex $w$ is called the central vertex. In the second
construction, the vertices $w_i$ are called the intermediate
vertices.
\begin{figure}[h]
  \begin{center}
    \includegraphics[height=3cm, width=12cm]{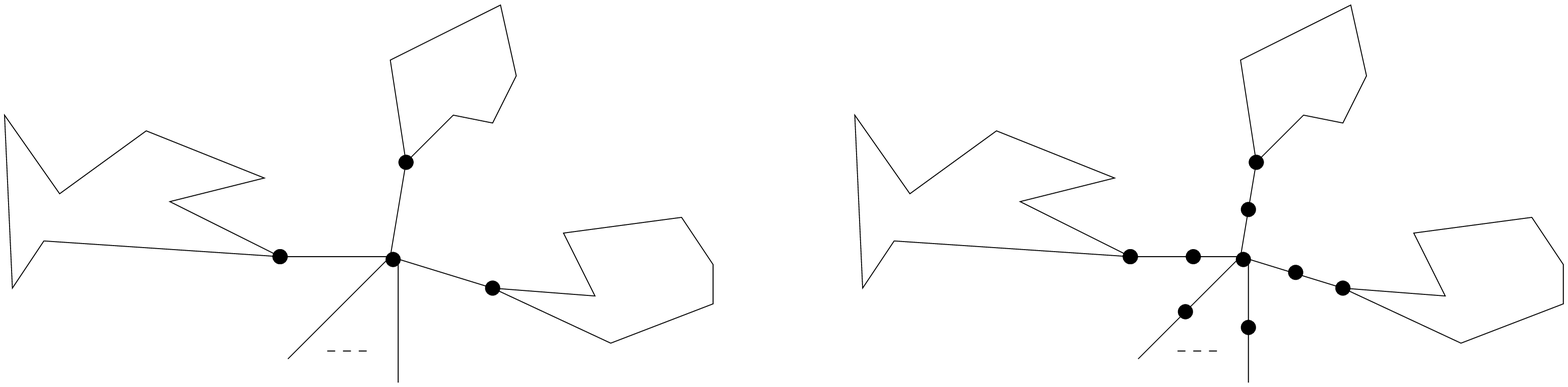}
            \begin{picture}(0,0)
                \put(-255,29){$w$}
                \put(-268,45){$v_i$}
                \put(-230,70){$G_i$}
                \put(-69,30){$w$}
                \put(-68,40){$w_i$}
                \put(-82,50){$v_i$}
                \put(-43,70){$G_i$}
                \put(-340,-15){Figure 21a: Simple building.}
                \put(-170,-15){Figure 21b: Complicated building.}
            \end{picture}
  \end{center}
\end{figure}
\begin{lemma}
Let $G_i^{v_i}$ be a brick of type $(n_i,m_i)$ for every $i\in
\{1, \dots ,n \}$. Let $S=S(G_1^{v_1},\dots ,G_k^{v_k})$ and
$C=C(G_1^{v_1},\dots ,G_k^{v_k})$. Then:

(1) $S^w$ is a brick of type $(\prod _{i=1}^{k}n_i-\prod
_{i=1}^{k}m_i, \prod _{i=1}^{k}n_i)$, where $w$ is the central
vertex of $S$.

(2) $S^{v_j}$ is a brick of type $(\prod _{i=1}^{k}n_i-\prod
_{i=1}^{k}m_i, m_j\prod _{j\neq i=1}^{k}n_i-\prod _{i=1}^{k}m_i)$.

(3) $C^w$ is a brick of type $(\prod _{i=1}^{k}(n_i-m_i)-\prod
_{i=1}^{k}n_i, \prod _{i=1}^{k}(n_i-m_i))$, where $w$ is the
central vertex of $C$.

(4) $C^{v_j}$ is a brick of type $(\prod _{i=1}^{k}(n_i-m_i)-\prod
_{i=1}^{k}n_i, -m_j\prod _{j\neq i=1}^{k}n_i)$.

(5) $C^{w_j}$ is a brick of type $(\prod _{i=1}^{k}(n_i-m_i)-\prod
_{i=1}^{k}n_i, n_j\prod _{j\neq i=1}^{k}(n_i-m_i)-\prod
_{i=1}^{k}n_i)$, where $w_j$ is any intermediate vertex of $C$.
\end{lemma}
\begin{proof}
We prove {\it (1)} and leave the other proofs to the reader:
\begin{eqnarray*}
f(S-w)&=&\prod _{i=1}^{k}f(G_i^{v_i})\ \ \ \ \ \ \ \ \ \
(multiplication)\\ &=&\prod _{i=1}^{k}n_i
\end{eqnarray*}
and
\begin{eqnarray*}
f(S)&=&f(S-w)-f(S-Nw)\ \ \ \ \ \ (recursion)\\ &=&\prod
_{i=1}^{k}n_i-\prod _{i=1}^{k}f(G_i-v_i)\ \ \ \ \ \ \ \
(multiplication)\\ &=&\prod _{i=1}^{k}n_i-\prod _{i=1}^{k}m_i.
\end{eqnarray*}
\qquad \qquad \qquad \qquad \qquad \qquad \qquad \qquad \qquad
\qquad \qquad \qquad \qquad \qquad \qquad \qquad \ \ \ \
\end{proof}

\vspace{0.3cm}

Let us come back to the example with $a_M=41$. Note that
$41=5\times 5 \times 2-3\times 3\times 1$. Hence $S=S(G_1^{v_1},
G_2^{v_2}, H^v)$ has $f(S)=41$ if $G_1^{v_1}$ and $G_2^{v_2}$ are
two copies of a brick of type $(5,3)$. Two different bricks of
type $(5,3)$ are shown in Figure~{22}. It follows from {\it (1)}
in above lemma that the planar graph $S$ convertible in BCG shown
in Figure 23 has $f(S)=41$. The corresponding diagram with
$a_M=41$ in its Kauffman bracket has $43$ crossings.

\begin{figure}[h]
  \begin{center}
    \includegraphics[height=1cm, width=12cm]{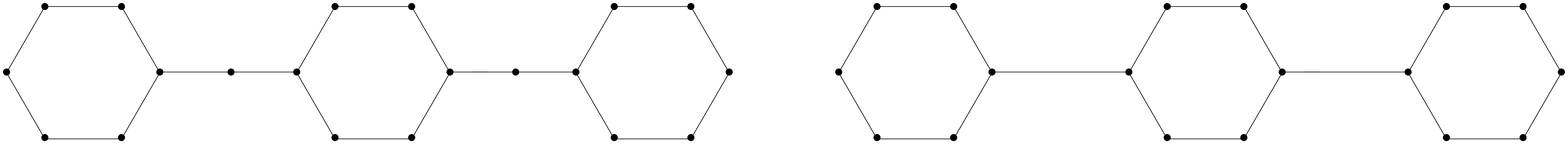}
    \caption{Two different bricks of type $(5,3)$.}
  \end{center}
\end{figure}
\begin{figure}[h]
  \begin{center}
    \includegraphics[height=2.5cm, width=7cm]{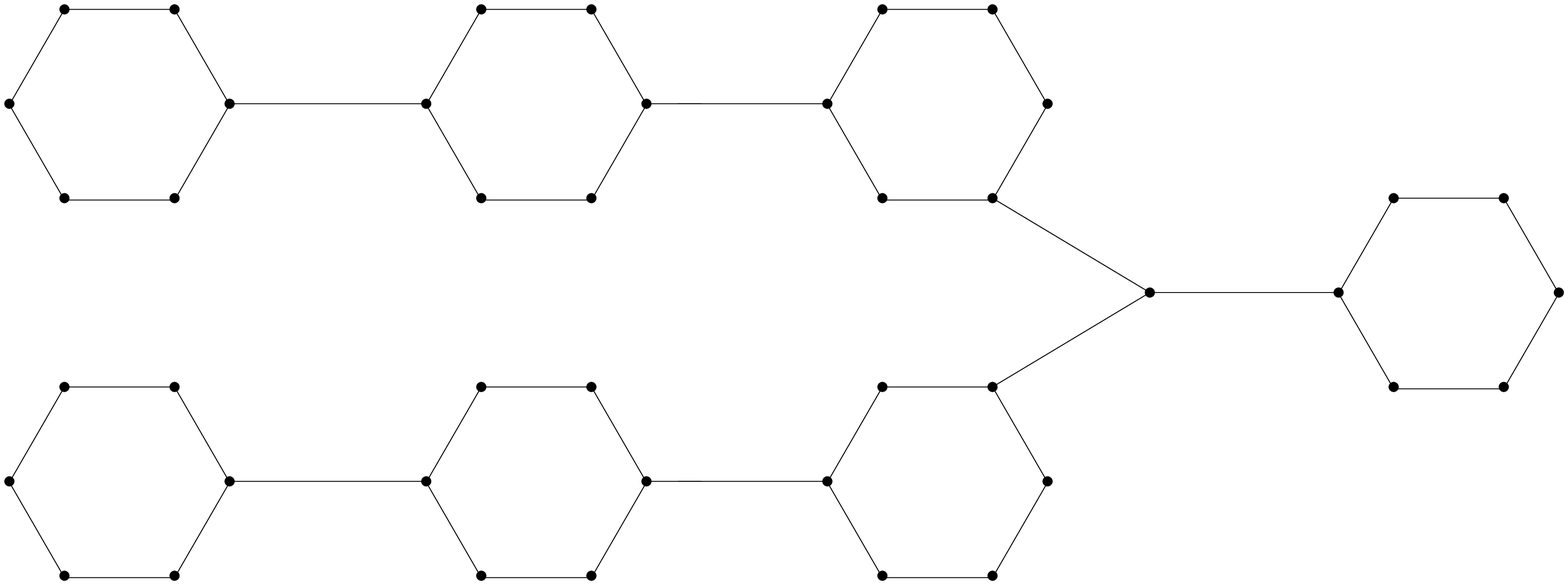}
    \caption{A graph $G$ with $f(G)=41$.}
  \end{center}
\end{figure}
It can be easily checked that any prime number less or equal than
$50$, different from $41$, can be realized as $f(S)$ where $S$ is
a simple building constructed with at most three hexagons (a brick
of type $(2,1)$).

\vspace{0.3cm}

Finally in this section we present a sequence $\{ F_r\} _{r\in
{\bf N}}$ of planar graphs, all of them convertible in BCG, such
that the corresponding sequence of integers $\{ f(F_r)\} _{r\in
{\bf N}}$ is the Fibbonachi sequence $2,3,5,8,13,21, \dots $. The
graph $F_r$ is shown in Figure 24. If we compare the link diagram
arising from $F_{r+1}$ and $F_r$, we note that while $a_M$
increases by $f(F_{r-1})$, the number of crossings of the diagram
increases only by $7$.
\begin{figure}[h]
  \begin{center}
    \includegraphics[height=1cm, width=6cm]{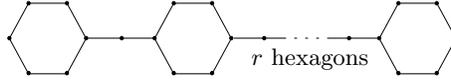}
    \caption{The graph $F_r$.}
        \begin{picture}(0,0)
                \put(7,38){{\small $r$ hexagons}}
        \end{picture}
  \end{center}
\end{figure}

\section{Non-states extreme coefficients of Jones polynomials.}

\vspace{0.3cm}

In this final section we will exhibit examples of prime knot
diagrams for which the extreme states coefficients $a_m$ and $a_M$
are zero and the next coefficients $a_{m+4}$ and $a_{M-4}$ take
arbitrary values. The idea is just to look at the more general
bipartite circle graph defined by $sD$, rather than the Lando's
graph $G_D$, and use a very simple trick for counting $a_{M-4}$ in
special cases.

We start with the natural enlargement of proposition given in the
second section. Recall that $\langle D \rangle$ denotes the
Kauffman bracket of an unoriented diagram $D$ with normalization
one.

\begin{proposition} \label{continuation}
\noindent (i) A state $s$ contribute to $a_{M-4}$ if and only if
either $s\in\Gamma _A=\{ s/|s|=|s_A|+b(s) \}$ or $s\in\Gamma
_A^1=\{ s/|s|=|s_A|+b(s)-2 \}$. The contribution of $s\in\Gamma
_A$ to $a_{M-4}$ is $(-1)^{|s_A|-1}(-1)^{b(s)}(|s_A|+b(s)-1)$. The
contribution of $s_1\in\Gamma _A^1$ to $a_{M-4}$ is
$(-1)^{|s_A|-1}(-1)^{b(s_1)}$.

\vspace{0.3cm}

\noindent (ii) $$a_{M-4}=(-1)^{|s_A|-1}[(|s_A|-1)\sum_{s\in \Gamma
_A}(-1)^{b(s)}+\sum_{s\in \Gamma _A}(-1)^{b(s)}b(s)+\sum _{s_1\in
\Gamma _A^1}(-1)^{b(s_1)}].$$

\vspace{0.3cm}

\noindent (iii) If $a_M=0$ then the whole contribution of $\Gamma
$ to $a_{M-4}$ is given by

$$(-1)^{|s_A|-1}\sum_{s\in \Gamma _A}(-1)^{b(s)}b(s).$$
\end{proposition}
\begin{proof}
A state $s$ contributes to $a_{M-4}$ if and only if $max(s)\geq
M-4$, hence if $s$ contributes to $a_{M-4}$ but is not in $\Gamma
_A$, we have that $max(s)=M-4$ by {\it (i)} in Proposition
\ref{well-known}. Now the statements follow from the fact that $\{
s_1/max(s_1)=M-4\} =\{ s_1/|s_1|=|s_A|+b(s_1)-2 \}$.
\end{proof}

\vspace{0.3cm}

Consider now the generalized bipartite circle graph $s_AD$
obtained from an unoriented diagram $D$ by $s_A$-smoothing.
Suppose that locally this has the very special aspect shown in
Figure 25, and suppose that only the A-chords $x_1, \dots ,x_k$
have their endpoints in the same component. The Lando's graph is
then given by $k$ parallel A-chords, hence by duplication we have
that $a_M=0$ in $\langle D \rangle $.
\begin{figure}[h]
  \begin{center}
    \includegraphics[height=3cm, width=12cm]{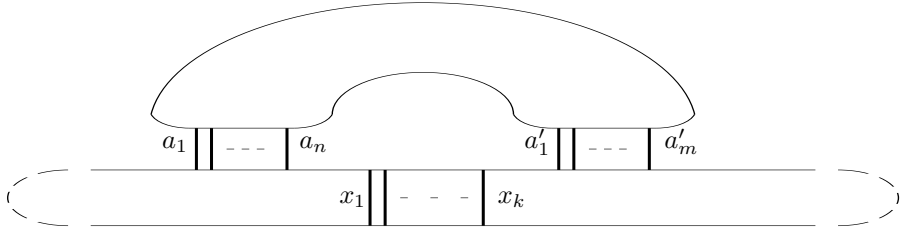}
    \caption{A very special $s_AD$.}
            \begin{picture}(0,0)
                \put(-110,64){$a_1$}
                \put(-58,64){$a_n$}
                \put(-43,43){$x_1$}
                \put(17,43){$x_k$}
                \put(27,64){$a_1'$}
                \put(80,64){$a_m'$}
            \end{picture}
  \end{center}
\end{figure}
From now on, we will identify a state $s$ with the set of the
$b(s)$ A-chords in $s_AD$ that correspond to the $b(s)$ crossings
of $D$ in which the label associated by $s$ is a B-chord. Let
${\cal X}=\{ x_1,\dots ,x_k \}$, ${\cal A}=\{ a_1,\dots ,a_n \}$
and ${\cal A}'=\{ a_1',\dots ,a_m' \}$ (see Figure 25).

\vspace{0.3cm}

The states in $\Gamma _A$ are the subsets of ${\cal X}$, hence by
{\it (iii)} in Proposition \ref{continuation} we have that the
contribution of $\Gamma _A$ to $a_{M-4}$ is given up to sign by
$\sum_{s\in \Gamma _A}(-1)^{b(s)}b(s)= \sum _{j=0}^{k} \left(
\begin{array}{c}
  k \\
  j
\end{array}
\right) (-1)^{j}j=0$.

We now analyze the contribution to $a_{M-4}$ of the states in
$\Gamma _A^1$ which contain at least one A-chord of ${\cal A}
\bigcup {\cal A}'$. Note first that these states are necessarily
contained in ${\cal X} \bigcup {\cal A} \bigcup {\cal A}'$.

Suppose first that in $s_1$ there are no A-chords of both sets
${\cal A}$ and ${\cal A}'$. Then the union of $s_1 \setminus {\cal
X}$ and any subset of ${\cal X}$ is a state which lies in $\Gamma
_A^1$, and all of their contributions to $a_{M-4}$ cancel.

By the contrary, suppose that $s_1$ has at least one A-chord of
${\cal A}$ and one of ${\cal A}'$. Then $s_1 \bigcap {\cal X}$ is
the empty set, and the contribution of all these states is
$(-1)^{|s_a|-1}$ in total.

\vspace{0.3cm}

As a consequence we have found a geometrical reason for some of
the results algebraically obtained in \cite{yo} about the Kauffman
bracket of pretzel link diagrams. Let $P$ be the pretzel link
diagram $P(a,b_1,\dots ,b_s)$ where $a\geq 2$, and $b_i\leq -2$
for any $i=1,\dots ,s$ with $s$ at least $2$ . Let $\alpha$ be the
cardinal of the set $\{ i / b_i=-2 \}$ and assume that $\alpha
\geq 1$. Then $P$ is a -adequate diagram, $a_M=0$ and
$a_{M-4}=\alpha$.

\vspace{0.3cm}

Finally, in order to get a prime knot diagram with extreme states
coefficient equal to zero and arbitrary values for $a_{m+4}$ and
$a_{M-4}$, we manipulate two pretzel link diagrams
$P(2,\underbrace{ -2,\dots ,-2}_{s})$ and $P(\underbrace{ 2,\dots
,2 }_{r-1},-2,2)$. We join these two diagrams in the way shown in
Figure 26, using two extra columns with $\alpha $ and $\beta $
crossings respectively, both $\alpha$ and $\beta $ great or equal
than two.
\begin{figure}[h]
  \begin{center}
    \includegraphics[height=3cm, width=12cm]{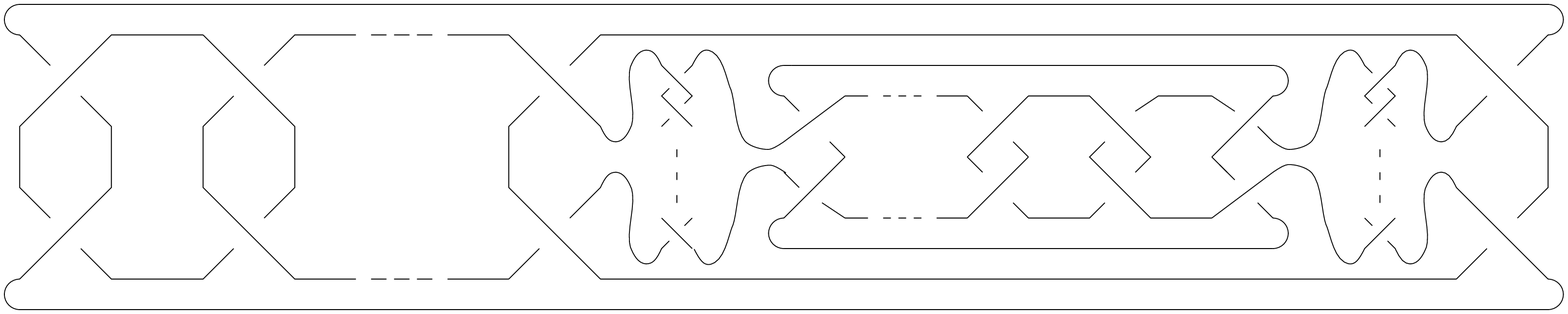}
    \caption{The link diagram $L(r,s;\alpha ,\beta )$.}
  \end{center}
\end{figure}
\noindent Note that this is the pretzel link diagram $P(2,
\underbrace{ -2, \dots , -2}_{s-2}, -\alpha , \underbrace{ 2,\dots
, 2}_{r-2}, -2, \beta )$ with two extra trivial knots placed in a
very special way, as shown in Figure 27.
\begin{figure}[h]
  \begin{center}
    \includegraphics[height=3cm, width=12cm]{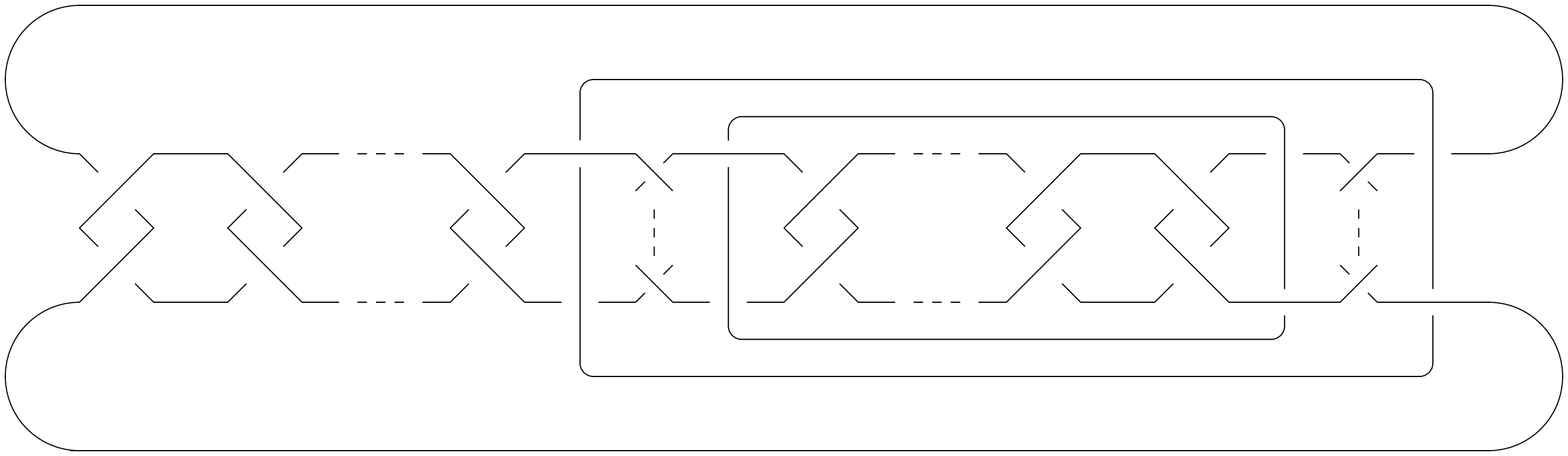}
    \caption{Other view of $L(r,s;\alpha ,\beta )$.}
  \end{center}
\end{figure}
\begin{theorem}
Let $L=L(r,s;\alpha ,\beta )$ be the link diagram shown in Figure
27. Then

$$\langle L \rangle =(-1)^{r+s+\beta-1}rA^{-4r-4s-\alpha -3\beta
+2}+\dots +(-1)^{r+s+\alpha -1}sA^{4r+4s+3\alpha +\beta -2}$$

\noindent In particular $span(\langle L \rangle )=8(r+s)+4(\alpha
+\beta )-4$. The highest degree of this polynomial is $M-4$, and
the lowest degree is $m+4$. Moreover, $L$ has
$r+s+2-\frac{1-(-1)^{\alpha }}{2}-\frac{1-(-1)^{\beta}}{2}$
components.
\end{theorem}
\begin{proof}Note first that $c(L)=2+2s+2+2r+\alpha +\beta$.

\vspace{0.3cm}

After $s_A$-smoothing (see Figure 28)we have $|s_A|=r+s+\alpha $
(hence $M=4r+4s+3\alpha +\beta +2$), and since the Lando's graph
is given by two parallel A-chords, we have $a_M=0$. The above
discussion applied to Figure 28 gives
$a_{M-4}=(-1)^{|s_A|-1}s=(-1)^{r+s+\alpha -1}s$.
\begin{figure}[h]
  \begin{center}
    \includegraphics[height=3cm, width=12cm]{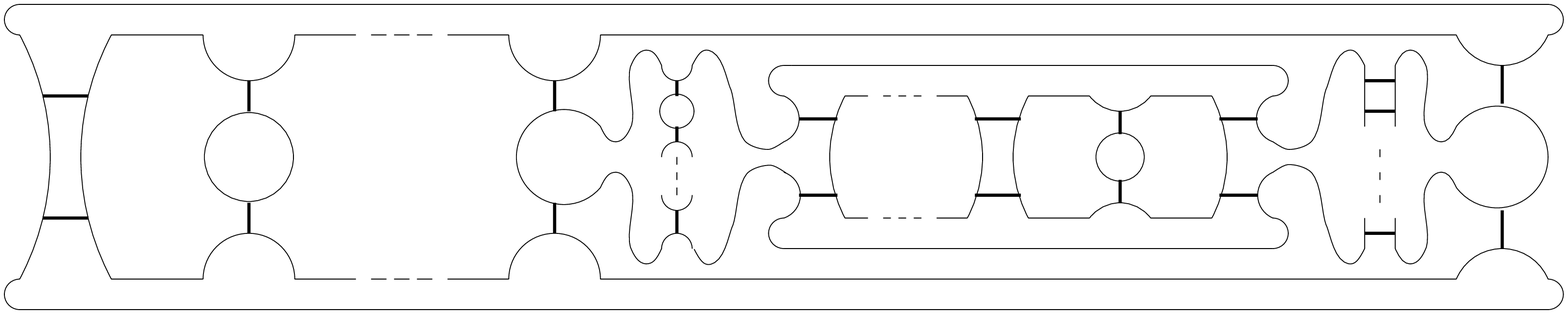}
    \caption{$L(r,s;\alpha ,\beta )$ afters $s_A$-smoothing.}
  \end{center}
\end{figure}
Analogously, after $s_B$-smoothing (see Figure 29) we have
$|s_B|=r+s+\beta $ (hence $m=-4r-4s-\alpha -3\beta -2$), and since
the Lando's graph is given by two parallel A-chords, we have
$a_m=0$. Finally, from the above discussion applied to Figure 29
we deduce that $a_{m+4}=(-1)^{|s_B|-1}r=(-1)^{r+s+\beta -1}r$.
\begin{figure}[h]
  \begin{center}
    \includegraphics[height=3cm, width=12cm]{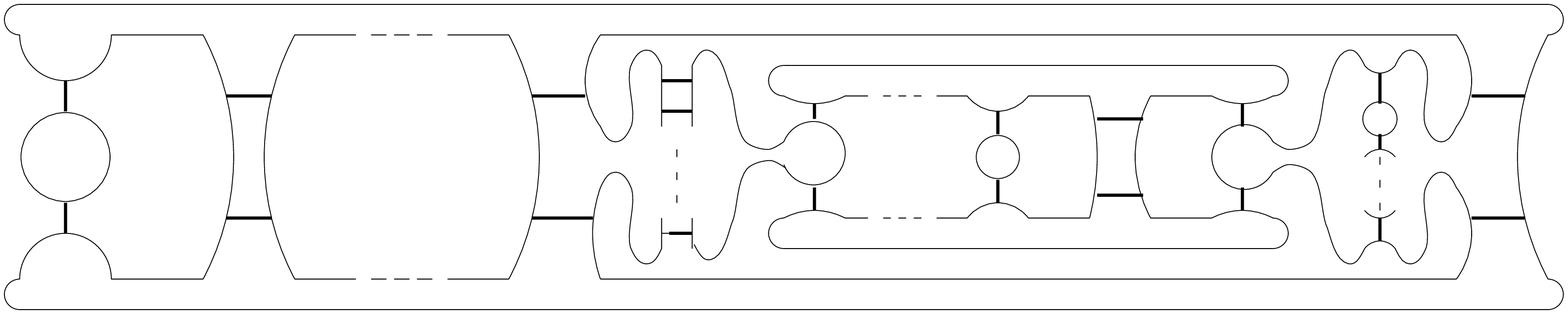}
    \caption{$L(r,s;\alpha ,\beta )$ after $s_B$-smoothing.}
  \end{center}
\end{figure}
\end{proof}

We now modify the link diagram $L(r,s;\alpha ,\beta )$ in order to
get a prime knot diagram $K(r,s;\alpha , \beta )$. This is shown
in Figure 30.
\begin{figure}[h]
  \begin{center}
    \includegraphics[height=3cm, width=12cm]{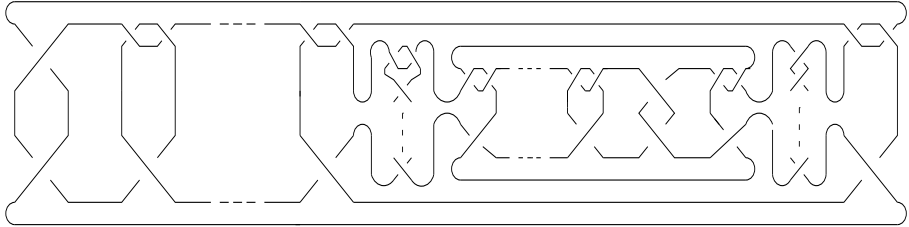}
    \caption{The diagram $K(r,s;\alpha ,\beta )$.}
  \end{center}
\end{figure}
\begin{theorem}
Let $K=K(r,s;\alpha ,\beta )$ be the link diagram shown in Figure
30. Then

$$\langle K \rangle =(-1)^{r+\beta}rA^{-5r-7s-\alpha -3\beta
-5}+\dots +(-1)^{s+\alpha -1}sA^{7r+5s+3\alpha +\beta +3}$$

\noindent In particular $span(\langle K \rangle )=12(r+s)+4(\alpha
+\beta )+8$. The highest degree of this polynomial is $M-4$, and
the lowest degree is $m+4$. Moreover, $K$ is a prime knot diagram.
\end{theorem}
\begin{proof}We first consider the modification $L'=L'(r,s;\alpha ,\beta )$ of
$L=L(r,s;\alpha ,\beta )$ shown in Figure 31.
\begin{figure}[h]
  \begin{center}
    \includegraphics[height=3cm, width=12cm]{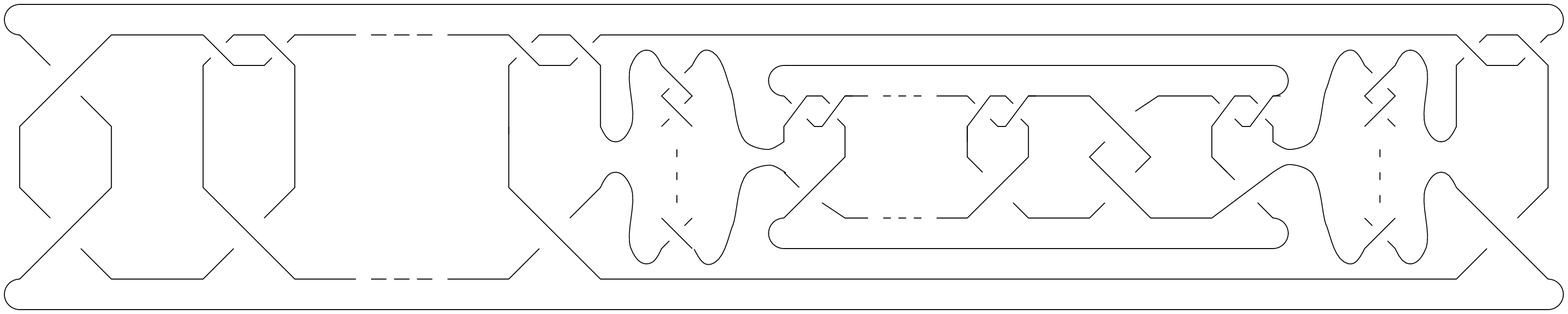}
    \caption{The link $L'=L'(r,s;\alpha ,\beta )$.}
  \end{center}
\end{figure}
Looking at the differences between $L$ and $L'$ after $s_A$ and
$s_B$-smoothings, we can prove that $\langle L' \rangle
=(-1)^{r+\beta -1}rA^{-5r-7s-\alpha -3\beta +2}+\dots
+(-1)^{s+\alpha -1}sA^{7r+5s+3\alpha +\beta -2}$ where the highest
and lowest degrees are $M-4$ and $m+4$ respectively.

\vspace{0.3cm}

But $L'$ is not a knot diagram when the parity of $\alpha $ and
$\beta $ is the same. For this reason we introduce a small change
in $L'$ in order to get the prime knot diagram $K$. This change is
shown in Figure 32.
\begin{figure}[h]
  \begin{center}
    \includegraphics[height=2cm, width=6cm]{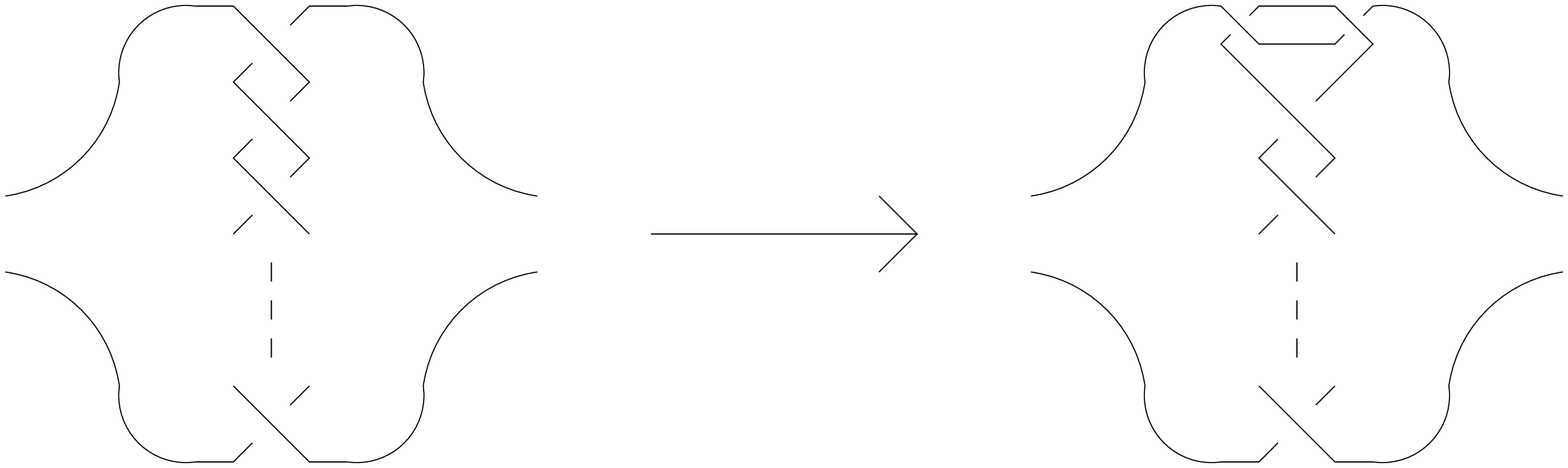}
    \caption{A small change in $L'$ gives $K$.}
            \begin{picture}(0,0)
                \put(-52,52){{\scriptsize $\alpha$}}
                \put(39,52){{\scriptsize $\alpha \!\!-\!\!1$}}
            \end{picture}
  \end{center}
\end{figure}
Note then that $c(K)=c(L')+1$ and $|s_A|$ does not change, hence
$M(K)=M(L')+1=7r+5s+3\alpha +\beta +3$. Respect to the
coefficients we have $a_M=0$ and
$a_{M-4}=(-1)^{|s_A|-1}s=(-1)^{s+\alpha -1}s$.

But after $s_B-$smoothing we have $|s_B|=|s_B(L')|+1$, hence
$m(K)=m(L')-1-2\times1=-5r-7s-\alpha -3\beta -5$. Respect to the
coefficients we have $a_m=0$ and
$a_{m+4}=(-1)^{|s_B|-1}r=(-1)^{r+\beta }r$.
\end{proof}

\vspace{0.3cm}

As we said in the introduction, we do not know if there is a nice
interpretation of $a_{M-4}$ in terms of graph theory, parallel to
that one in which $a_M$ is given in terms of the graph $G_D$ as
described in the second section. Apart from the examples, this
last section can be seen as a partial interpretation of $a_{M-4}$
in terms of graph theory. More important, it remains to be
answered the following question: how arbitrary can be the extreme
coefficients of the Jones polynomial when the value for the spread
is previously fixed?

\end{document}